\documentclass[a4paper,11pt]{amsart}
\usepackage{amsthm,amsmath,amssymb}
\usepackage[OT1]{fontenc}
\usepackage{graphicx,wrapfig}
\usepackage{hyperref}
\usepackage[utf8]{inputenc}
\usepackage{enumitem}

\newtheorem*{theorem}{Theorem}

\let\leq\leqslant
\let\geq\geqslant

\newcommand{\first}{\ensuremath{1^\text{st}}}
\newcommand{\second}{\ensuremath{2^\text{nd}}}

\makeatletter
\let\old@setaddresses\@setaddresses
\def\@setaddresses{\bgroup\parindent 0pt\let\scshape\relax\old@setaddresses\egroup}
\makeatother

\begin{document}

\title{A Note on Concurrent Graph Sharing Games}

\author[S.~Chaplick]{Steven Chaplick}
\author[P.~Micek]{Piotr Micek}
\author[T.~Ueckerdt]{Torsten Ueckerdt}
\author[V.~Wiechert]{Veit Wiechert}

\address[S.~Chaplick, V.~Wiechert]{Institut f\"ur Mathematik, Technische Universit\"at Berlin, Germany}
\email{\{chaplick,wiechert\}@math.tu-berlin.de}
\thanks{S.~Chaplick is supported by the ESF research project EUROGIGA GraDR.}

\address[P.~Micek]{Theoretical Computer Science Department, Faculty of Mathematics and Computer Science, Jagiellonian University, Poland}
\email{piotr.micek@tcs.uj.edu.pl}
\thanks{P. Micek is supported by the Mobility Plus program from The Polish Ministry of Science and Higher Education.}

\address[T.~Ueckerdt]{Department of Mathematics, Karlsruhe Institute of Technology, Germany}
\email{torsten.ueckerdt@kit.edu}


\thanks{V.\ Wiechert is supported by the Deutsche Forschungsgemeinschaft within the research training group `Methods for Discrete Structures' (GRK 1408).} 

\begin{abstract}
 In the concurrent graph sharing game, two players, called \first{} and \second{}, share the vertices of a connected graph with positive vertex-weights summing up to $1$ as follows. 
 The game begins with \first{} taking any vertex. 
 In each proceeding round, the player with the smaller sum of collected weights so far chooses a non-taken vertex adjacent to a vertex which has been taken, i.e., the set of all taken vertices remains connected and one new vertex is taken in every round. 
 (It is assumed that no two subsets of vertices have the same sum of weights.)
 One can imagine the players consume their taken vertex over a time proportional to its weight, before choosing a next vertex.
 In this note we show that \first{} has a strategy to guarantee vertices of weight at least $1/3$ regardless of the graph and how it is weighted.
 This is best-possible already when the graph is a cycle.
 Moreover, if the graph is a tree \first{} can guarantee vertices of weight at least $1/2$, which is clearly best-possible.
\end{abstract}

\maketitle


Imagine a pizza, sliced as usually into triangular pieces, not necessarily of the same size, and two players alternatingly taking slices in such a way that every slice, except the first one, is adjacent to a slice that was taken earlier.
What is the fraction of the total size of the pizza that the first player can guarantee to get at least, independently of the number of slices and their sizes (weights)?
This problem, the so-called \emph{Pizza Problem}, posed by Peter Winkler was resolved in~\cite{CKM+10,KMU11} and it turns out that \first{} can always guarantee to get at least $4/9$ of the entire pizza and that this is best-possible.
Considering a pizza to be a cycle with weights on its vertices, one can find work on similar games for trees~\cite{MiW11,SeS12} and subdivision-free graphs~\cite{GMW14+}.

The \emph{concurrent graph sharing game} is a variant of the Pizza Problem introduced by Gao in~\cite{Rao} (as \emph{the Pizza Race Game} and its generalizations).
As before, a vertex-weighted graph is shared by \first{} and \second{} taking one vertex at a time in such a way that the set of all taken vertices remains connected.
The game begins with \first{} taking any vertex. 
However in each proceeding round, the player with the smaller sum of collected weights so far picks the next non-taken vertex.
(We assume for now that no two subsets of vertices have the same sum of weights and discuss the situation without this assumption at the end of the paper.)
One can imagine the players consume their taken vertex over a time proportional to its weight, before choosing a next vertex.

For convenience, assume that the weights of all the vertices in the graph sum up to $1$.
In~\cite{Rao} the author claims that for every weighted cycle \first{} can guarantee to take vertices of total weight at least $2/5$.
However, his proof has a flaw and cannot be fixed.
In fact, he starts by introducing a new vertex with vanishingly small weight between any two adjacent vertices of the given cycle, and claims that these vertices are irrelevant for the analysis of the game.
But this is true only as long as \first{} does not start with such a vertex.
(When \first{} starts with an original vertex Gao's argument for a $2/5$ lower bound seems to be correct).
Indeed, we show here that the maximum total weight that \first{} can guarantee on every cycle is $1/3$.
In fact, our lower bound argument works for every graph, i.e., \first{} can always guarantee to take vertices of total weight at least $1/3$.

Secondly, Gao asks whether \first{} can guarantee any positive fraction of the total weight if the game is played on a tree.
We show here with an easy strategy stealing argument that, playing on trees, \first{} can always guarantee to take vertices of total weight at least $1/2$, which is clearly best-possible.

An \emph{instance} of the concurrent graph sharing game is a pair $(G,w)$ of a graph $G = (V,E)$ and positive real vertex weights $w:V \to (0,1]$ with $\sum_{v \in V} w(v) = 1$.
For a subset $A \subseteq V$ of vertices we denote $w(A) = \sum_{a \in A} w(a)$.
For a vertex $a \in V$, let $F_a$ and $S_a$ be the subsets of vertices that \first{} and \second{} take when \first{} starts with $a$, and from then on both players play optimally subject to maximizing $w(F_a)$ and $w(S_a)$, respectively. 
Thus, $F_a\sqcup S_a = V$ for all $a\in V$.
The \emph{value of an instance $(G,w)$} is the maximum total weight $v(G,w)$ of vertices that \first{} can guarantee to take in this instance.
In particular, $v(G,w) = \max_{a \in V(G)} w(F_a)$.

\begin{theorem}\label{thm:main}
 For the concurrent graph sharing game we have
 \[
  \inf_{(G,w)} v(G,w) = 1/3 \quad \text{ and } \quad \inf_{(G,w), \text{ $G$ is a tree}} v(G,w) = 1/2.
 \]
\end{theorem}

Recall that in the above Theorem we consider only instances $(G,w)$ in which no two disjoint subsets of vertices have the same weight.
However as explained below, we use this hypothesis only for a strategy stealing argument proving the lower bound of $1/2$ for trees.
We remark that a similar strategy stealing is part of Gao's proof.

\begin{proof}
 To prove $\inf_{(G,w)} v(G,w) \geq 1/3$, let $(G,w)$ be any instance of the concurrent graph sharing game.
 If there is a vertex $a \in V(G)$ with $w(a) \geq 1/3$, then clearly $v(G,w) \geq w(F_a) \geq 1/3$.
 On the other hand, if $w(a) < 1/3$ for all $a \in V(G)$, then at the moment \first{} can take no further vertex (because all vertices are already taken), \second{}'s current vertex has weight less than $1/3$. 
 So for every $a \in V(G)$ we have $w(S_a) - w(F_a) < 1/3$.
 Together with $w(F_a) + w(S_a) = 1$ this implies that $v(G,w) \geq w(F_a) > 1/3$.
 
 \medskip

 Next we shall prove $\inf_{(G,w)}v(G,w) \leq 1/3$ by providing for every $\varepsilon > 0$ an instance $(G,w_\varepsilon)$ with $v(G,w_\varepsilon) \leq 1/3 + \varepsilon$.
 Consider the cycle $G$ consisting of seven vertices $a,b,c,d,e,f,g$ in this cyclic order and corresponding vertex-weights $M, M+15, 17, 7, 12, M+26, 18$, where $M = M(\varepsilon) \gg 95$ is large enough.
 This instance $(G,w_\varepsilon)$ is depicted in Figure~\ref{fig:worst-pizza-ever} in form of a pizza.
 It contains three pieces with weight at least $M$, which we call \emph{heavy}.
 For \second{} to get at least two heavy pieces (and therefore roughly $2/3$ of the entire pizza) he moves according to the table on the right and then takes the last heavy piece when it is his turn again.
 E.g., when \first{} starts with $d$, \second{} takes $e$, and if \first{} continues with $c$, \second{} takes $f$, and then \second{} is guaranteed to get the last heavy piece (either $a$ or $b$ in this case).
 
 \begin{figure}[htb]
  \centering
  \includegraphics{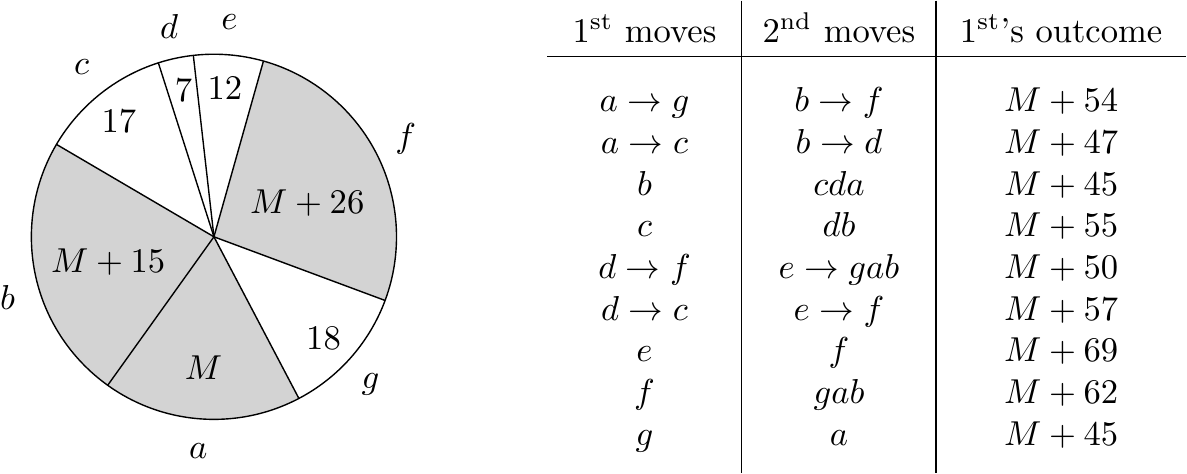}
  \caption{A pizza (vertex-weighted cycle) with total weight $3M + 95$ in which \first{} cannot guarantee to get more than $M + 69$, for any $M$ large enough. Divide each weight by $3M+95$ to get a total sum of weights equal to $1$.}
  \label{fig:worst-pizza-ever}
 \end{figure}

 We now consider instances in which the graph is a tree. 
 To prove that $\inf_{(G,w), \text{ $G$ tree}} \geq 1/2$ let $(G,w)$ be any instance where $G$ is a tree.  
 If $|V(G)| = 1$, then clearly $v(G,w) = 1$.
 Otherwise, for each vertex $a \in V(G)$, let $b(a) \in Y_a$ be the first vertex \second{} takes when \first{} starts with $a$.
 As $|E(G)| < |V(G)|$ there exists an edge $aa'$ such that $b(a) = a'$ and $b(a') = a$.
 Consider the games in which \first{} starts with $a$ and $a'$, respectively, and both players play optimally.
 In the former game \first{} starts with $a$ and \second{} answers with $a'$, while in the latter it is the other way around.
 In particular, from that moment on both games are identical, but the roles of \first{} and \second{} are switched.
 It follows that $w(X_a) = w(Y_{a'})$ and with $w(Y_{a'}) = 1 - w(X_{a'})$ we conclude $v(G,w) \geq \max\{w(X_a), w(X_{a'})\} \geq 1/2$.
 
 \medskip
 
 To see that $\inf_{(G,w), \text{ $G$ is a tree}} \leq 1/2$, consider for every $\varepsilon > 0$ a tree consisting of a single edge $ab$ with $w(a) = (1 - \varepsilon)/2$ and $w(b) = (1 + \varepsilon)/2$.
\end{proof}

\begin{figure}[htb]
 \centering
 \includegraphics{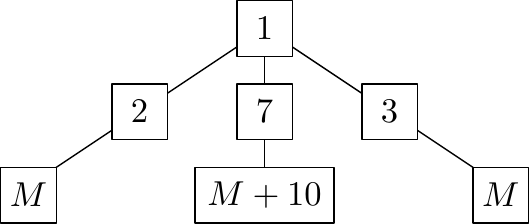}
 \caption{A vertex-weighted tree with total weight $3M + 23$ in which \first{} cannot guarantee to get more than $M + 20$, for any $M$ large enough. Divide each weight by $3M+23$ to get a total sum of weights equal to $1$.}
 \label{fig:worst-tree-ever}
\end{figure}

We remark that the assumption that no two subsets of vertices have the same sum of weights is crucial for the strategy stealing argument we used for the trees.
Indeed, if both players may finish their current vertex at the \emph{same} time and in such cases, say, \first{} always takes the next vertex, then the best \first{} can guarantee on any tree is $1/3$, instead of $1/2$.
For the upper bound see Figure~\ref{fig:worst-tree-ever}.
On the other hand, note that the general lower bound of $1/3$ remains valid, no matter how those ``ties'' are broken.

\medskip

Lastly, we remark that in the original Pizza Problem most of the effort was to find a strategy for \first{} to get at least $4/9$ of the pizza.
The concurrent graph sharing games considered in this paper turned out to be somewhat simpler in the analysis.
Indeed, here the strategy for \first{} to get $1/3$ of the total weight in any graph is obvious and the difficulty was to believe that this is best-possible already for cycles.
Tight examples were found by working out the sequence of moves in an optimal strategy for \second{}, which led to a system of linear constraints for the vertex-weights whose optimization gave a best-possible scenario for \second{}.

\bibliographystyle{plain}
\bibliography{pizza}
\end{document}